\documentclass[12pt]{amsart}

\usepackage{amsmath}   \usepackage{amscd}
\usepackage{graphics}  \usepackage{latexsym}
\usepackage{amsfonts}  \input{psfig.sty}

\theoremstyle{plain}              
\newtheorem{Thm}{Theorem} \newtheorem{Lem}[Thm]{Lemma}      \newtheorem{Cor}[Thm]{Corollary}
                         \theoremstyle{definition}
                                
\theoremstyle{remark}     

\def\Z{\mathbb Z}

\def\cp{\hbox{${\mathbb C} P^2$}}
\def\cpb{\hbox{$\overline{{\mathbb C}P^2}$}}

\def\s2x{\hbox{$S^2 \times S^2$}}

\def\sqr#1#2{{\vcenter{\hrule height.#2pt    		
\hbox{\vrule width.#2pt height#1pt \kern#1pt       		
\vrule width.#2pt}\hrule height.#2pt}}}	

\def\square{\mathchoice\sqr67\sqr67\sqr{2.1}6\sqr{1.5}6}

\def\qed{~\hfill$\square$}

\begin{document}
\title[]{An involution acting nontrivially on Heegard-Floer homology}
\author{Selman Akbulut and Selahi Durusoy} 
\thanks{First named author is partially supported by NSF grant DMS 9971440}
\address{Department of Mathematics \\ Michigan State University \\ MI, 48824}
\email{akbulut@math.msu.edu \and durusoyd@math.msu.edu}
\date{\today}
\begin{abstract}We show that a certain involution on a homology sphere $\Sigma$  induces a nontrivial homomorphism on its  Heegard-Floer homology groups (recently defined by Ozsv\'ath and Szab\'o). We discuss  application of this to constructing exotic smooth structures on $4$-manifolds. 
\end{abstract}

\maketitle
\vspace{-0.2in}
\setcounter{section}{-1}

%
%

\section{Introduction}Let $W$ be the smooth contractible $4$-manifold consisting of a single $1$- and $2$- handle attached as in Figure 1 (one of the Mazur manifolds).
\begin{figure}[ht]  \begin{center}    
 \includegraphics{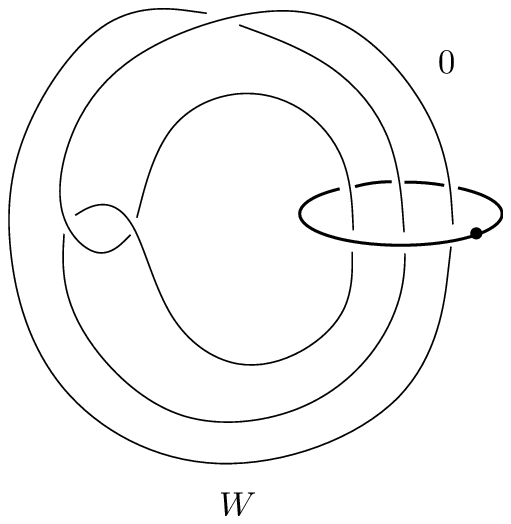}   \caption{} \label{mazur1}    \end{center}  \end{figure}

\noindent The boundary is a homology sphere $\Sigma= \partial W$. Let $f: \Sigma \to \Sigma$ be the obvious involution obtained by first surgering $S^{1}\times B^{3}$ to $B^{2}\times S^{2}$ in the interior of $W$, then surgering the other imbedded $B^{2}\times S^{2}$ back to  $S^{1}\times B^{3}$ (i.e. replacing the dots in \mbox{Figure \ref{mazur1}).} This involution interchanges the small linking loops of the two circles in Figure \ref{mazur1}.  
It is known that this involution acts nontrivially on the Donaldson-Floer homology of $\Sigma$, (c.f. \cite{a1}, \cite{s}). 
 In this paper we will show that the involution $f$ also acts nontrivially on the Heegard-Floer homology group $HF^+(\Sigma)$ of $\Sigma$, which was defined in \cite{os1}. We will deduce the proof from the construction of \cite{a2} and the $4$-manifold invariants of \cite{os3}.  This gives a possible way  of changing smooth structures of $4$-manifolds that contain $W$ (such codimension zero contractible submanifolds are usually called ``corks''\!, e.g.  \cite{k}, \cite{m}, \cite{s}, and they can be made Stein \cite{am}).

\begin{figure}[ht]  \begin{center} 
 \includegraphics{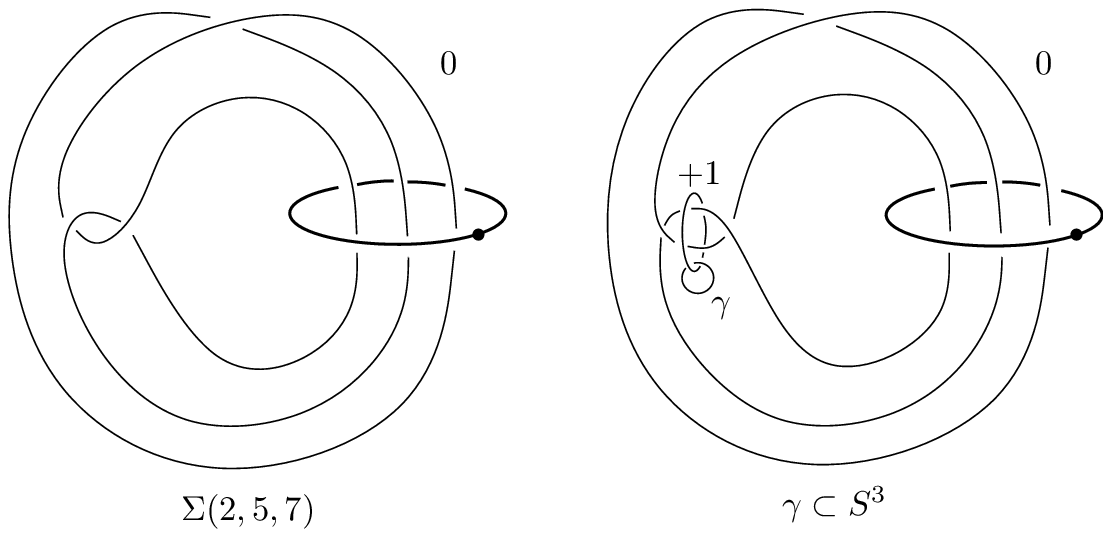}   \caption{} \label{mazur2}    \end{center}  \end{figure}

\begin{figure}[ht]  \begin{center}     \includegraphics{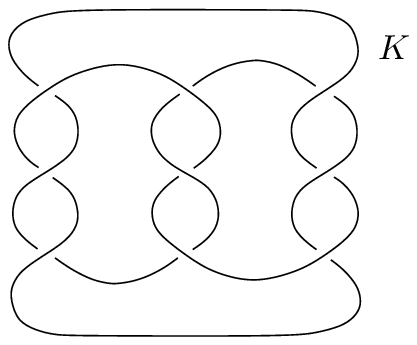}   \caption{} \label{pretzel_knot}    \end{center}  \end{figure}

By \cite{ak}, the boundary of the $4$-manifold in the first picture of the Figure \ref{mazur2} is the Brieskorn homology sphere $\Sigma(2,5,7)$, which we can relate to $\Sigma$ as follows: The second picture of Figure \ref{mazur2} is obtained by attaching a $+1$ framed $2$-handle to $W$, hence the boundary is  $S^{3}$ (since the $0$-framed $2$- handle  slides over the $+1$ framed handle, and cancels the $1$-handle). So, in this picture, $\gamma$ is just a loop in $S^{3}$,
which is drawn nonstandardly.  It is easily seen that in the standard picture of $S^{3}$, the loop $\gamma$ corresponds to the $(-3,3,-3)$ pretzel knot $K \subset  S^{3}$ as shown in Figure \ref{pretzel_knot} (e.g. \cite{ak}).

 Note that doing $r+1$ surgery to  $\gamma$ corresponds to doing $r$ surgery to  $K$ in $S^{3}$, which we will denote by $S^{3}_{r}(K)$. Hence $S^{3}_{1}(K) =\Sigma(2,5,7)$ and $S^{3}_{-1}(K) =\Sigma $. We will use this identifications in the Heegard-Floer homology calculations of the next section.

%
%

\section{Calculating the Heegaard-Floer Homology $HF^+( \Sigma )$}

We first calculate the Heegaard-Floer homology of $\Sigma(2,5,7)$
using techniques of \cite{os5},                  
from this and the surgery exact sequence of \cite{os2} 
we will deduce the Heegaard-Floer homology of $\Sigma$.

From \cite{ak} 
we know that $\Sigma(2,5,7)$ is the boundary of the negative definite plumbing of disk bundles
over $2$-spheres described by the following plumbing graph:

%
\begin{figure}[ht]
 \centerline{\psfig{figure=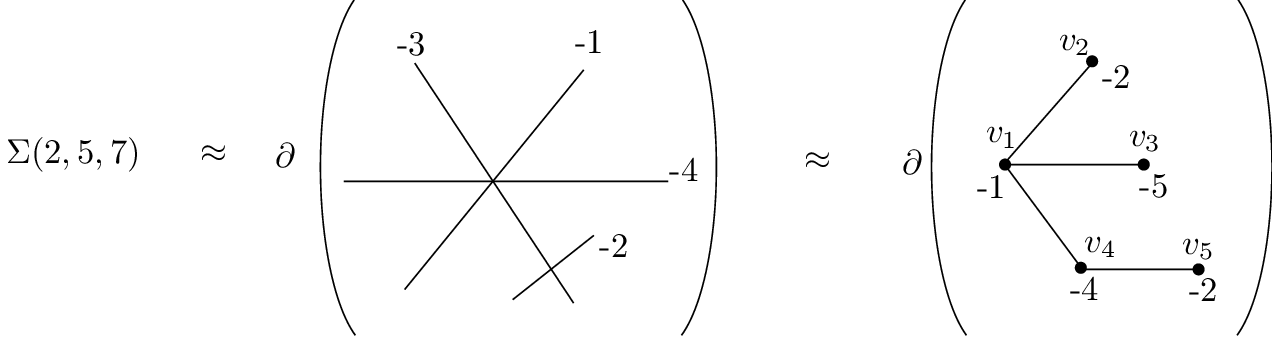}}
 \caption{}
 \label{sigma257_graph}
\end{figure}

This graph is negative definite and has only one bad vertex $v_1$ in the sense of \cite{os5}, 
i.e., only $[v_1] \cdot [v_1] > -d(v_1)$ where $d(v)$ counts edges containing $v$.

In  the notation of \cite{os5} we write $\Sigma(2,5,7) = \partial X(G) = Y(G)$, and  we compute 
$$HF^+(-Y(G)) \cong H^+(G)$$
where
$H^+(G)$ is a group which can be calculated  by the  algorithm of  \cite{os5} for negative definite (possibly disconnected) trees $G$ with at most one bad vertex. Here we  follow this algorithm (for the terminology, we refer to
\cite{os5}): 

First we need to find characteristic vectors $K$ satisfying
\begin{equation}
 [v]\cdot[v] + 2 \leq \langle K,[v] \rangle \leq -[v]\cdot[v].
 \label{cond13} 
\end{equation}

We remark that there are only finitely many such vectors since $G$ is a definite graph, in particular there are $80$ characteristic vectors satisfying (\ref{cond13}) for the given graph. Proposition 3.2 in \cite{os5} characterizes
a spanning set for $K^+(G)$ which is the dual point of view for $H^+(G)$. Among the 80 characteristic vectors satisfying
(\ref{cond13}) we look for those vectors K which carry full paths ending at vectors $L$ for which $-L$ satisfies (\ref{cond13}).
There are precisely three such vectors K and we list them and also the full paths for convenience (paths are obtained by consecutively adding $2PD(v_i)$ for the numbers $i$ listed):

\begin{tabular}{rcll}
$K_1$ & = & $12v_1 + 6v_2 + 3v_3 + 4v_4 + 2v_5$ & path: 1,2,1  \\
$K_2$ & = & $-8v_1 - 4v_2 - v_3  - 2v_4 - 2v_5$ & path: 1,2,1,5,4,1,2,1,3,1,2,1,4,1,2,1,5  \\
$K_3$ & = & $-16v_1- 8v_2 - 3v_3 - 4v_4 -2v_5$  & path: 1,2,1
\end{tabular}

\vspace{.1in}

All other equivalence classes in $K^+(G)$ can be obtained using the $U$-action from these vectors. 
We must check which are distinct. From now on, we will write the vectors as $5$-tuples with entries
$ \langle K, [v_i] \rangle $, hence $K_i$ will be denoted as
$$(1,0,-3,-2,0), \;(1,0,-3,-2,2), \text{ and } (1,0,-1,-2,0).$$

The following diagram shows how $U\otimes K_1$, $U \otimes K_2$ and $U \otimes K_3$ are equivalent to $(-3,4,1,2,0)$,
and hence to each other:

\begin{figure}[ht]
 \centerline{\psfig{figure=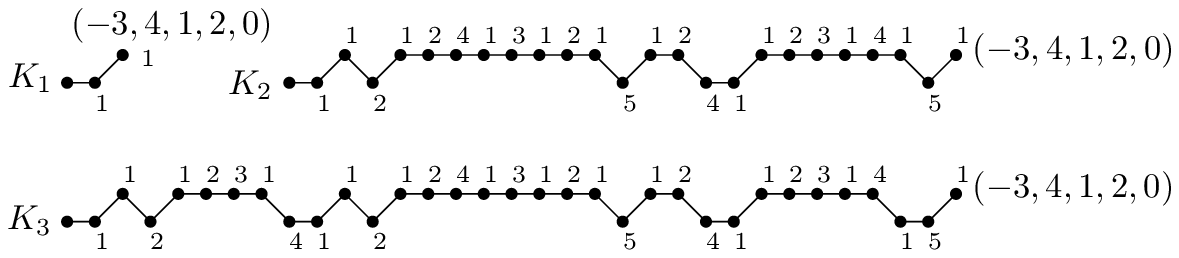}}
 \caption{}
\end{figure}

In the above figure, vertices correspond to characteristic vectors and an edge ending with a number $i$ means
$2PD(v_i)$ was added to previous vector to reach the new one. Edges with positive slope mean tensoring with $U$.
\vspace{.05in}
\noindent Hence first diagram represents the vectors
$(1,  0, -3, -2,  0),  (-1,  2, -1,  0,  0),  (-3,  4,  1,  2,  0)$ which tells us that
$$(1,  0, -3, -2,  0) \sim (-1,  2, -1,  0,  0) \text{ and }
 U \otimes (-1,  2, -1,  0,  0) \sim (-3,  4,  1,  2,  0)$$

\noindent Hence each $K_i$ lies on the same grading level and we can check this by calculating
$$\frac{K\cdot K + |G|}{4}$$ which gives zero for each $K_i$.
Therefore we conclude:
$HF^+(-Y(G)) \cong T_0^+ \oplus \mathbb{Z}_{(0)} \oplus \mathbb{Z}_{(0)}$
By Prop 7.11 of \cite{os3} we can convert this to the Heegard-Floer homology of $Y(G)$:
$$HF^+_{*}( Y(G)) = HF_{-}^{-*-2}(-Y(G)) = T_0^+ \oplus \mathbb{Z}_{(-1)} \oplus \mathbb{Z}_{(-1)} $$ 
To do this, we first use  the following long exact sequence  to compute $HF^{-}(-Y(G))$:

$$...\to HF^-(-Y(G)) \longrightarrow  HF^\infty (-Y(G)) \longrightarrow  HF^+ (-Y(G))\to ...$$

\noindent then we adjust gradings by changing the signs and substracting $-2$. Note that here we adapted  the convention of  \cite{os5} by denoting $T_s^+ =  \bigoplus _{k=0}^{\infty} \Z_{(s + 2k)}$. So we have:
$$HF^+(\Sigma (2,5,7))= T_0^+ \oplus \mathbb{Z}_{(-1)} \oplus \mathbb{Z}_{(-1)}$$

\vspace{.1in}

{\Lem  $HF^+(\Sigma) \cong T_0^+ \oplus \mathbb{Z}_{(0)} \oplus \mathbb{Z}_{(0)}$.}

\proof Let us consider the Heegard-Floer homology long exact sequence of \cite{os2}  for the $(-3,3,-3)$ pretzel knot $K\subset S^3$ :
$$...\to HF^+(S^3) \longrightarrow  HF^+(S^{3}_{0}(K)) \longrightarrow  HF^+ (S^{3}_{1}(K))\to ...$$
 Recall that  the middle two arrows in the sequence decrease  the gradings  by $1/2$.
Since $ S^{3}_{1}(K) =\Sigma(2,5,7)$ we know $2$ out of  $3$  terms of the exact sequence, from this we can compute the third term  $HF^+ (S^{3}_{0})(K)=T^{+}_{(1/2)}\oplus T^{+}_{(-1/2)}\oplus  \Z_{(-1/2)} \oplus \Z_{(-1/2)} $,    as shown in the diagram:

$$\begin{array}{ccccclc}

&  ...  \to HF^+(S^3) & \longrightarrow  & HF^+(S^{3}_{0}(K)) & \longrightarrow  &  HF^+ (\Sigma(2,5,7))   \to..&  \\
        &   :                   &                                    & :                        &                               &                                    & \\
        &       \Z_{(4)}   &   \vector (4,-1){30}    & \Z _{(9/2)}      & \vector (4,-1){30}     & : & \\
       &                       &                                       &   \Z _{(7/2)}      &                       & \Z_{(4)}& \\
       &       \Z_{(2)}   &      \vector (4,-1){30}    &  \Z _{(5/2)}       &                          \vector (4,-1){30}        &   & \\
       &                       &                                       & \Z _{(3/2)}        &                  &  \Z_{(2)} & \\
       &      \Z_{(0)}    &    \stackrel{g}{ \vector (4,-1){20} }   &     \Z _{(1/2)}         & \vector (4,-1){30}&  & \\
      &                       &        \mbox{nonzero}  &  \Z _{(-1/2)} \oplus \Z _{(-1/2)} \oplus \Z _{(-1/2)}   & \vector (4,-1){20} &\Z_{(0)}& \\
      &                       &                                    &                       &  &   \Z_{(-1)} \oplus  \Z_{(-1)}& \\
      \end{array}$$
      
      \vspace{.1in}

\noindent Now by plugging in this value of $HF^+ (S^{3}_{0} (K))$ in the Heegard-Floer homology exact sequence  of the knot $K$ in $S^3$ below (this time involving $-1$, $0$, and no surgeries),  we calculate  the Heegard -Floer homology of $ \Sigma = S^{3}_{-1}(K)$ to be 
$$HF^+(\Sigma) \cong T_0^+ \oplus \mathbb{Z}_{(0)} \oplus \mathbb{Z}_{(0)}$$
$$\begin{array}{rrccllll}

        & ...  \to HF^+(S^{3}_{-1}(K)) & \longrightarrow  & HF^+(S^{3}_{0}(K)) & \longrightarrow  &  HF^+ (S^3) & \to..   \\
        &   :                   &                                    & :                        &                               &                                    & \\
        &       \Z_{(4)}   &   \vector (4,-1){30}    & \Z _{(9/2)}      & \vector (4,-1){30}     & : & \\
       &                       &                                       &   \Z _{(7/2)}      &                       & \Z_{(4)}& \\
       &       \Z_{(2)}   &      \vector (4,-1){30}    &  \Z _{(5/2)}       &                          \vector (4,-1){30}        &   & \\
       &                       &                                       & \Z _{(3/2)}        &                  &  \Z_{(2)} & \\
       &  \Z_{(0)}\oplus \Z_{(0)} \oplus \Z_{(0)}   &     \vector (4,-1){20}    &     \Z _{(1/2)}         & \vector (4,-1){30}&  & \\
      &                       &            &  \Z _{(-1/2)} \oplus \Z _{(-1/2)} \oplus \Z _{(-1/2)}   & & \Z_{(0)}  \;\; \stackrel{h} {\vector (4,-1){20}}&
      & \\
      &&&&&\;\;\;\;\;\;\;\;\mbox{zero}&
       \end{array}$$
      
      \vspace{.1in}
      
In these calculations we used the fact that all the  groups are equipped  with $U$-action (an action of $\Z [U]$ which lowers the grading by $-2$), and the maps are equivariant with respect  to this action.  

Also we used  two additional facts, which were explained to us by P. Ozsv\'ath:   The lower left map $ g $ of  the first diagram is nonzero, and the  lower right map $h$ of the second diagram is zero.  Explanation: By using Theorem 1.4 of \cite{os6} and Theorem 6.1 of \cite{os7}  we can calculate $HF^+(S^{3}_{0}(K))$ independently, and by the exact sequence and the fact that $g$ has to lower the degree by $1/2$ implies that $g$ must be nonzero. The map $h$ is zero because it is obtained by summing over the maps induced by a cobordism $(Z,s)$ with $Spin^{c}$ structures $s$ (extending the one on the boundary).  $Z$ has signature $-1$ and  there are two such $Spin^c$ structures which cancel each other. \qed

%
%
\section{ Heegard-Floer homology action on  $W$}

  Let $W_{0}$ be  the punctured $W$  (i.e. $W- B^{4}$).  We want to compute the map induced by the cobordism $W_{0}$  \cite{os2} 
 $$F^+_{W_{0}}: HF^+(\Sigma) \to HF^+(S^3).$$
 
\begin{figure}[ht]  \begin{center}  
\includegraphics{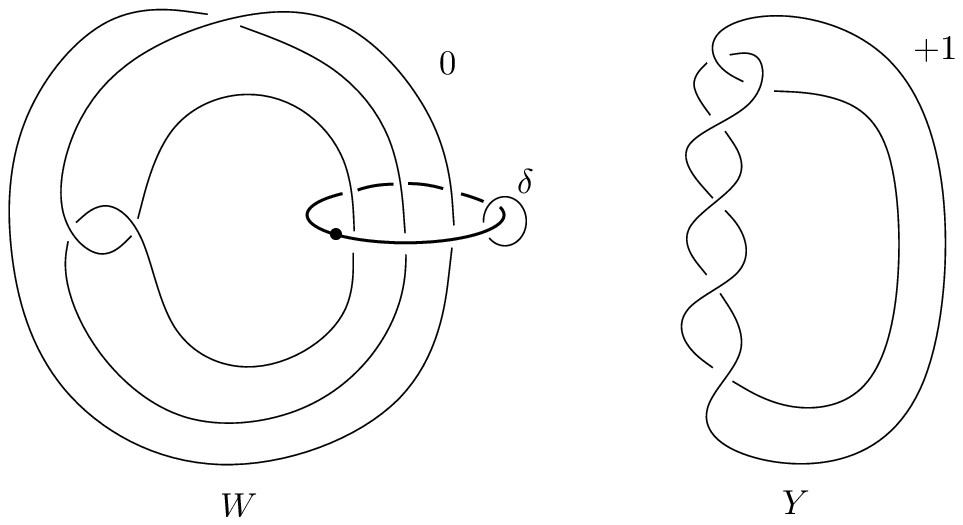}   \caption{} \label{surgery_on_delta}    \end{center}
   \end{figure}

$W$ consists of $S^1\times B^3$ with a $2$-handle attached. By turning  $W$ upside down we see  $W_{0}$ is a union of  cobordisms $W_{1}$ from $\Sigma $ to $S^1 \times S^2 $,  and  $W_{2}$  from $S^1\times S^2$ to $S^3$ (a 3-handle). Hence we have
 the decomposition $F^+_{W_{0}}=F^+_{W_{2}} \circ F^+_{W_{1}}$. The map
 $F^+_{W_{2}} : HF^+(S^1\times S^2) \to HF^+(S^3)$ is well understood, it is obtained by the projection 
  $$T^{+}_{ 1/2} \oplus T^{+}_{-1/2}\to T^{+}_{-1/2}\stackrel {\cong}{ \to}T^{+}_{0}$$
 The second map is $U$-equivariant,  in particular it maps the lowest degree generator $\Z_{(-1/2)}$ of  $T^{+}_{-1/2}$ to the lowest degree generator $Z_{(0)}$ of $T^{+}_{0}$. So it is sufficient to understand the map
 $F^+_{W_{1}}: HF^+(\Sigma) \to HF^+(S^1\times S^2) $. Notice $W_1$ is obtained by attaching
a $2$-handle to $\Sigma = \partial W$ along the loop $\delta$ with $0$-framing (see Figure \ref{surgery_on_delta}).

Denote  the manifold obtained from $\Sigma$, doing $k$-surgery to $\delta$, by $\Sigma_{k}(\delta)$. Note that $\Sigma_{0}(\delta)=S^1\times S^2$  and $\Sigma_{1}(\delta)$ is the second manifold of Figure \ref{surgery_on_delta}, which we will call $Y$ to simpilfy the notation. 
We can write the corresponding Heegard-Floer homology exact sequence  of $\delta \subset \Sigma$:

\begin{equation} ...\to HF^+(\Sigma) \longrightarrow  HF^+(\Sigma_{0}(\delta)) \longrightarrow  
 HF^+ (Y)\to ...
\label{delta_surgery_exact_seq} 
\end{equation} 
 
\begin{figure}[ht]  \begin{center} 
 \includegraphics{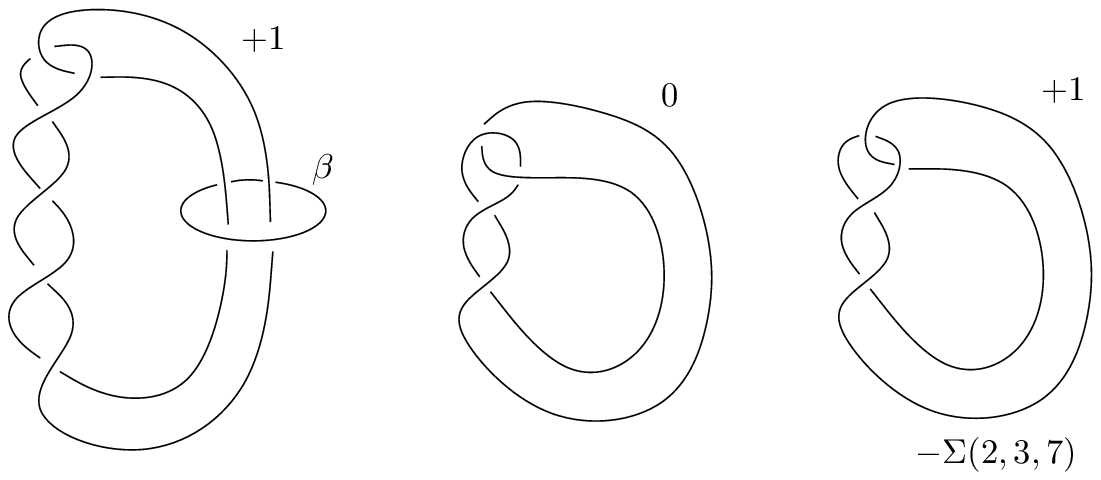}   \caption{} \label{surgery_on_beta}   \end{center}  \end{figure}

%
%
$F^+_{W_{1}}$  is just the first map in the exact  sequence, to compute it we need to know the third term in the  sequence. To calculate $HF^+ (Y)$  we use  the Heegard-Floer exact sequence of the  knot $\beta \subset Y$ of Figure \ref{surgery_on_beta}.
$$...\to HF^+(Y) \longrightarrow  HF^+(Y_{0}(\beta)) \longrightarrow  
 HF^+ (Y_1 (\beta)) \to ...$$
It is easy to check that $Y_{0}(\beta)$ is the manifold obtained by $0$-surgery to the figure eight knot, and $Y_{1}(\beta)$ is $-\Sigma (2,3,7)$. In \cite{os4}, \cite{os5} homologies of these  were computed : $HF^+(Y_{0}(\beta))=T^+_{1/2}\oplus T^+_{-1/2}\oplus \Z_{(-1/2)} $ and 
$HF^+(-\Sigma(2,3,7))=T^{+}_{0}\oplus \Z_{(0)}$. By plugging these to the last exact sequence we compute
$$HF^+(Y)= T^+_{0}\oplus Z_{(0)}\oplus\Z_{(0)} $$
By  plugging the value $HF^+(Y)$  in the previous exact  sequence (\ref{delta_surgery_exact_seq})   we see that the map  
$ F^+_{W_{1}}: HF^+(\Sigma) \to HF^+(S^1\times S^2) $ is given by the obvious projection and the inclusion  (where the middle map is an $U$-equivariant  isomorphism) $$T^{+}_{0}\oplus\Z_{(0)}\oplus \Z_{(0)} \to T^{+}_{0} \stackrel {\cong}{\to} T^{+}_{-1/2} \to T^{+}_{-1/2}\oplus T^{+}_{1/2}$$
Hence the map $ F^+_{W_{0}}: HF^+(\Sigma) \to HF^+(S^3) $ is given by the obvious projection
\begin{equation}
T^{+}_{0}\oplus\Z_{(0)}\oplus \Z_{(0)} \to T^{+}_{0}
 \end{equation}

%
%

\section{ Nontriviality of the involution $f$}

Let $Z$ be the $K3$ surface constructed in \cite{a1}. It was shown that $Z \, \# \, \cpb $ decomposes as the union  of two codimension zero submanifolds glued along their common boundries $$Z \, \# \, \cpb =N\cup W$$
where $W$ is the Mazur manifold of Figure 1, furthermore  by identifying $N$ and $W$  via the involution
$f:\Sigma\to \Sigma$ (instead of the identity) gives the decomposed manifold $$N\cup_{f}W=3(\cp) \, \# \, 20 \cpb$$

Now let us remember the $4$-manifold invariants of \cite{os3}: Let $(X,s) $ be a smooth closed 4-manifold with a $Spin^c$ structure  $s$  such that it is decomposed  as a union of two codimension zero submanifolds glued along their boundaries $X=X_{1}\cup_{\partial} X_{2}$ with $b_{2}^{+}(X_{i}) >0 $ , $i=1,2$. Call $Y=\partial X_{i}$. By puncturing each $X_{i}$ in the interior  we get two cobordisms: $X'_{1}$ from $S^3$ to $Y$,  and $X'_{2}$ from $Y$ to $S^3$.

\begin{figure}[ht]  \begin{center} \includegraphics{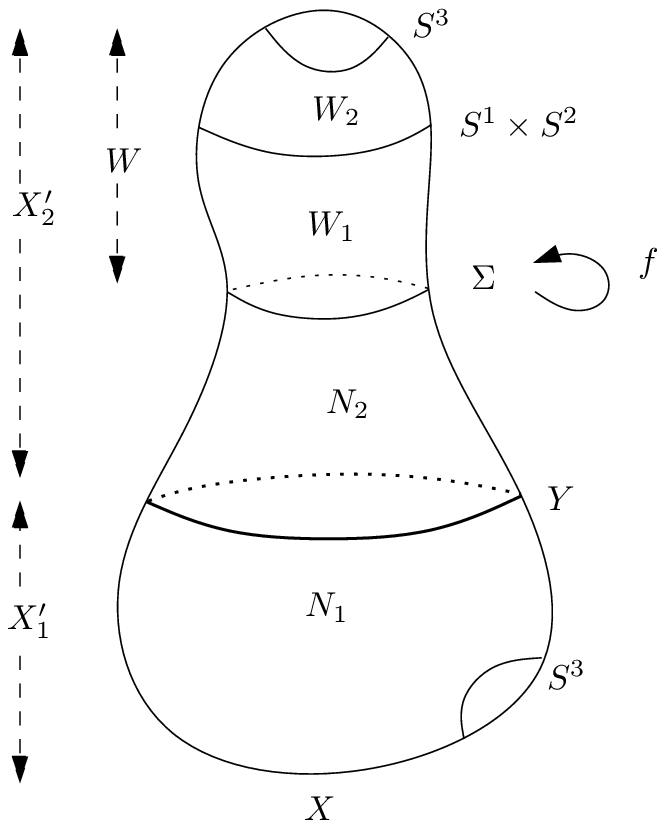}   \caption{} \label{}    \end{center}  \end{figure}

\noindent Since $b_{2}^+(X'_i)>0$ maps $F^{\infty}_{X'_{1}}$ and  
 $F^{\infty}_{X'_{2}}$ induced by cobordism are zero (upper right and lower left vertical arrows in the diagram), hence the following commuting diagram 
 
 $$\begin{array}{cccccccc}
... &\longrightarrow  &HF^+(S^3)   & \stackrel{\delta}{\longrightarrow} &HF^-(S^3)    &\longrightarrow     & HF^{\infty}(S^3) & \longrightarrow...\\
 &        &\big{\downarrow} &        \vspace{.05in}                       &   \text{\phantom{F\,X}} \big{\Downarrow}  F^{-}_{X_{1}'} &         & \text{\phantom{zero}} \big{\downarrow}  \mbox{\small{zero}} &\\
...\longrightarrow    HF^{\infty}(Y)    & \longrightarrow &HF^+(Y) &  \Longrightarrow                                      & HF^{-}(Y)     & \longrightarrow     &HF^{\infty}(Y) & \longrightarrow...   \\
  \text{\phantom{zero}} \big{\downarrow} \mbox{\small{zero}}  &    &  \text{\phantom{F\,X}} \big{\Downarrow}  F^{+}_{X_{2}'} &       \vspace{.05in}            & \big{\downarrow} &           & \big{\downarrow}  &  \\
...\longrightarrow  HF^{\infty} (S^3)&\longrightarrow & HF^+(S^3) & \longrightarrow                                 & HF^{-}(S^3) & \longrightarrow      &HF^{\infty}(S^3) & \longrightarrow ...  
\end{array}$$

\noindent induces a well defined map:
$F^{mix}_{(X,s)}: HF^{-}(S^3)\to HF^{+}(S^3)$. Then the Ozsv\'ath-Szab\'o   invariant $\Phi _{(X,s)}\in \Z$  is the degree of this map on $\Z_{(-2)}$ generator (in particular the invariant doesn't depend on the splitting of $X$). By   \cite{os4}, \cite{os8} we know that  this invariant  is zero on $3(\cp) \, \# \, 20 \cpb  $ while it is non-zero on $K3 \, \# \, \cpb$ (actually to get this conclusion we also need  to show that either the homotopy $K3$ of \cite{a2} is diffemorphic to actual $K3$, or directly apply the technique of \cite{os8} to show that it has nonzero invariant - either way works.). Also, since this invariant is independent of splittings, we can further split  
$X = Z \, \# \, \cpb $ as shown in Figure 8, and decompose $$ F^{+}_{X'_{2}}= F^{+}_{W_{0}}\circ F^{+}_{N_{2}} $$

 We can view the manifold $N\cup_{f}W$ as being obtained from  $N\cup W$ by cutting it along $\Sigma$ and sticking in the cobordism manifold $H=\Sigma\times [0,1/2]\cup_{f}\Sigma\times [1/2,1]$.
Also the  involution $f$ induces a map  $f^{*}:HF^+(\Sigma) \to HF^+(\Sigma) $ via the cobordism $H$,
($ f^*=F^+_{H} $ in \cite{os3}'s notation). So the 
nonzero invariant of  $Z \, \# \, \cpb $ is obtained by computing the degree of the  map: 
$$\Z_{(-2)} \stackrel {\lambda}{\to} T_{0}^+\oplus \Z_{(0)} \oplus \Z_{(0)}\to T_{0}^+ $$
(where $\lambda$ is some map and the second map is the projection) whereas the zero  invariant of 
$3(\cp) \, \# \, 20 \cpb$ is obtained by computing the degree of the map: 
$$\Z_{(-2)} \stackrel {\lambda}{\to} T_{0}^+\oplus \Z_{(0)} \oplus \Z_{(0)}
\stackrel{f^{*}}{ \to }T_{0}^+\oplus \Z_{(0)} \oplus \Z_{(0)} \to T_{0}^+ $$

Note that since $H$ is diffeomorphic to the product $\Sigma \times [0,1]$ the $U$-equvariant map $f^{*}$ should look like the identity map. At first glance this appears to contradict the above result. The explanation is that,  the two ends of $H$ are different copies of $\Sigma$, only when we fix the both ends by a reference copy of $\Sigma$, then $f^{*}$ becomes an isomorphism permuting two of the zero degree generators of 
$T_{0}^+\oplus \Z_{(0)} \oplus \Z_{(0)}$. More specifically,  $f^{*}$  permutes the zero degree generator of $T_{0}^+$ with one of the generators of the remaining  $\Z_{(0)} \oplus \Z_{(0)} $. So, actually $f^{*}$ is only a $U$-equivariant isomorphism, which happens to differ from the identity isomorphism. Hence we proved:
{\Thm $f: \Sigma \to \Sigma$ induces a nontrivial involution $f^{*}:HF^{+}(\Sigma )\to HF^{+}(\Sigma  )$}.

%
%

If $(M,s)$ is any smooth cobordism from $S^3$ to $\Sigma$ with $Spin^c$-structure $s$,  and  $b_{2}^{+}(M)>1$, then  Ozsv\'ath-Szab\'o procedure gives a map  $F^{mix}_{(M,s)}:HF^{-}(S^3)\to HF^{+}(\Sigma)$. Previous theorem gives a potential way of detecting exotic smoothings of $4$-manifolds.
 
{\Cor Let $Q=N\cup _{\partial} W$ be a smooth closed $4$-manifold with $b_{2}^{+}>1$ (the union is taken along the common boundary $\Sigma$). Let $N_{0}$ be the cobordism from $S^3$ to $\Sigma$  obtained from $N$ by removing a copy of  $B^4$ from its interior. Then  the manifold $Q'= N\cup_{f} W$ is a fake copy of $Q$, provided that the image of the map $F^{mix}_{(N_{0},s) }HF^{-} (S^3) \to HF^{+}(\Sigma)$ lies in $T^{+}_{0}$ for some $s$.} 

\vspace{.06in}

\noindent {\it Acknowledgements:} We would like to thank P.Ozsv\'ath for illuminating discussions on Heegard-Floer homology, and R.Kirby  for encouragement, and MSRI for providing stimulating environment where this research is completed.


\vspace{-.03in}

\begin{thebibliography}{99999}

\bibitem[A1]{a1}S. Akbulut, {\em An Involution permuting Floer Homology},
Turkish J. Math. 18 (1994) 16--22

\bibitem[A2]{a2}S. Akbulut, {\em A Fake Compact Contractible $4$-Manifold},
J. Differential Geom. 33 (1991) \mbox{335--356}

\bibitem[AK]{ak}S. Akbulut and R.Kirby, {\em Mazur Manifolds}, 
Michigan Math. J. 26 (1979) 259--284

\bibitem[AM]{am} S. Akbulut and R. Matveyev, {\em A Convex Decomposition Theorem for $4$-manifolds}, 
Internat. Math. Res. Notices No. 7 (1998), 371--381. arXiv:math.GT/0010166

\bibitem[K]{k} R. Kirby, {\em Akbulut's Corks and $h$-cobordisms of Smooth, Simply Connected $4$-manifolds},
\mbox{Turkish} J. Math. 20 (1996), 85--93

\bibitem[M]{m} R. Matveyev, {\em A decomposition of smooth simply-connected $h$-cobordant $4$-manifolds},
J. Diff. Geom. 44 (1996), 571--582


\bibitem[OS1]{os1} P. Ozsv\'ath and Z. Szab\'o,
{\em Holomorphic Disks and Topological Invariants for Closed Three-Manifolds},
arXiv:math.SG/0101206 v4

\bibitem[OS2]{os2} P. Ozsv\'ath and Z. Szab\'o,
{\em Holomorphic Disks and Three-Manifold Invariants: Properties and Applications},
arXiv:math.SG/0105202 v1

\bibitem[OS3]{os3} P. Ozsv\'ath and Z. Szab\'o,
{\em Holomorphic Triangles and Invariants for Smooth Four-Manifolds},
arXiv:math.SG/0110169 v2

\bibitem[OS4]{os4} P. Ozsv\'ath and Z. Szab\'o,
{\em Absolutely Graded Floer Homologies and Intersection Forms for Four-Manifolds with Boundary},
arXiv:math.SG/0110170 v2

\bibitem[OS5]{os5} P. Ozsv\'ath and Z. Szab\'o, {\em On the Floer Homology of Plumbed Three-Manifolds},
Geom. Topol. 7 (2003) 185--224

\bibitem[OS6]{os6} P. Ozsv\'ath and Z. Szab\'o, {\em Heegard-Floer homology and Alternating Knots},
Geom. Topol. 7 (2003) 225--254

\bibitem[OS7]{os7} P. Ozsv\'ath and Z. Szab\'o, {\em Knot Floer Homology, Genus Bounds, and Mutation},
arXiv:math.GT/0303225 v1

\bibitem[OS8]{os8} P. Ozsv\'ath and Z. Szab\'o, {\em Holomorphic Triangle Invariants and the Topology of Symplectic Four-Manifolds} arXiv:math.SG/0201049 v1 (2002)


\bibitem[S]{s} N. Saveliev, {\em A Note on Akbulut Corks},
Math. Res. Lett. 10 (2003) 777--785

\end{thebibliography}
\end{document}